%% file: main.tex
\documentclass[11pt]{article}

\usepackage{geometry}
 \usepackage[belowskip=-10pt,aboveskip=0pt]{caption}
 
\usepackage[utf8]{inputenc}
\usepackage{amsmath}
\usepackage{xcolor}
\usepackage{graphicx}
\usepackage{float}
\usepackage{subcaption}
\usepackage[ruled]{algorithm2e}
\usepackage[english]{babel}
\usepackage{blindtext}
\usepackage{hyperref}
\usepackage{url}
\usepackage{indentfirst}
\usepackage{amsbsy}

\definecolor{forgreen}{RGB}{34,139,34}
\newcommand{\ket}[1]{\ensuremath{|#1\rangle}}

\DeclareMathOperator*{\argmin}{argmin}

\title{Adaptive Strategies For Efficient Model Reduction In High-Dimensional Inverse Problems}
\author{Andrei Mukhin, Aleksey Khlyupin}
\date{December 2018}

\begin{document}
\maketitle


\begin{abstract}
\input{abstract.tex}
\end{abstract}


\maketitle

\section{Introduction}

    \input{Introduction/introduction.tex}

\section{Gradient-based history matching}
    \input{Section2/subsection2_1.tex}
    \input{Section2/subsection2_2.tex}
    \input{Section2/subsection2_3.tex}
    
\section{Adaptive parametrization using AS-PCA}
    \input{Section3/subsection3_1.tex}

    \input{Section3/subsection3_2.tex}

    \input{Section3/subsection3_3.tex}
\section{Results}
\label{sec::results}
    \subsection{Synthetic problem description}
        \input{Results/TPD.tex}

\section{Discussion}
    \input{conclusion.tex}

\bibliography{biblio} 
\bibliographystyle{abbrv}
\end{document}

%% file: abstract.tex
This work explores a novel approach for adaptive, differentiable parametrization of large-scale non-stationary random fields. Coupled with any gradient-based algorithm, the method can be applied to variety of optimization problems, including history matching. The developed technique is based on principal component analysis (PCA), but, in contrast to other PCA-based methods, allows to amend parametrization process regarding objective function behaviour.

To define an efficient basis update, Adaptive Strategies PCA (AS-PCA) uses gradients of objective function with respect to model parameters. Gradients are already available in gradient-based optimization process combined with adjoint method. Optimality, correctness and low computational cost of the new parametrization procedure is guaranteed by theoretical derivation of basis with stationary perturbation theory from quantum mechanics. Three modifications of method are presented. One of them not only improves quality of parametrization, but also extends applicability of method for uncertainty quantification.

The AS-PCA is then applied to synthetic history matching. Forward problem is represented by non-linear parabolic equation. Through simple cases, we demonstrate an improve in matching performance in comparison with standard linear PCA for misfit minimization and overall field consistency. 

%% file: Introduction/introduction.tex
Inverse problems appear in many areas of comparative research, where the problem of defining inner structure of uncertain object having a set of its lifecycle observations is considered. Although correct solution for a lot of practical inverse problems requires an efficient parametrization algorithm, our specific interest here is in the field of history matching problem. The purpose of this procedure is to generate detailed reservoir description that is consistent with prior information and match production data to within some tolerance. This is usually done using two types of data, namely static and dynamic. Static data is mostly constant over time, e.g. geological concept of formation, well logs and petrophysical data and is commonly given as prior information. Dynamic data is time-dependent and represents properties that change during production process, e.g. pressure and flow rates, flow responses. Initial reservoir model is constructed with static data and then has to be adjusted to correctly match dynamic data. This is done by iterative change of model parameters (e.g. facies or permeability) untill model flow simulation is in well agreement with production data. A relatively recent review on history matching problem can be found in \cite{oliver2011recent}.

A common approach is to perform history matching either within optimization formulation, or as data assimilation problem. For the latter, ensemble-based methods, such as ensemble Kalman filters (EnKF) \cite{aanonsen2009ensemble}, recently have gained popularity. Such methods are non-invasive regarding forward flow simulator and provide multiple matched models. These features simplify process of uncertainty quantification and allow to work with black-box simulator. However, among limitations of ensemble-based methods is problem of ensemble collapse. This leads to need of large number of realizations within ensemble and, in hence, high computational costs. Also, these methods are incapable to represent models with complex geological structures.

In optimization context, the history matching problem is often addressed by stochastic methods such as genetic algorithm \cite{romero2000modified}, particle swarm optimization \cite{hajizadeh2011towards}, evolutionary algorithms \cite{schulze2001optimization} and others. Although these methods perform a global search and allow to use forward simulator as black-box, their use for complex models has a relatively high computational cost and can be applicable only via decent computational clusters.

A gradient-based optimization approach is in the main focus of this paper. These methods are invasive w.r.t the forward simulator and provide a local search, but they are sufficiently faster than methods described above. Adjoint-based techniques are usually applied in history matching procedure in this context, since they provide neccesary gradients at cost of one additional simulation. These techniques are investigated for partial differential equations (PDE) based systems, e.g. with applications in closed-loop reservoir management \cite{sarma2006efficient} and, recently, for integro-differential equation (IDE) systems with a memory effect \cite{kadyrova2018application}. 

Since history matching often has to be performed on the real large-scale fields, parametrization techniques are very useful as they substantially reduce the number of parameters that have to be determined. Also, it allows for maintaining geological consistency of result, to some extent. Approaches like discrete cosine transform (DCT) \cite{jafarpour2008history} or discrete wavelet transform (DWT) \cite{sahni2005multiresolution} allow to represent model in term of relatively few parameters, but their performance on complex models is not suitable since information about prior covariance of geological model is not considered in basis construction process. Another choice is to use PCA-based algorithms. Classic linear PCA, often referred to as Karhunen-Loeve expansion in this context, has been succesfully applied to history matching problem \cite{sarma2006efficient}. Due to the fact that linear PCA considers only covariance matrix of model, it preserves only two-point geostatistics, and number of approaches were developed to use PCA for complex non-Gaussian fields. Among them are kernel PCA (kPCA) \cite{sarma2008kernel}, optimization-based PCA (O-PCA) \cite{vo2014new}, regularized kPCA (R-kPCA) \cite{vo2016regularized}, convolutional neural network PCA (CNN-PCA) \cite{liu2018deep}.

One of the missing key in existing methods is the lack of ability to adaptively include information about objective function structure into parametrization process. Sensitivity of data is considered in a number of parametrization techniques \cite{vogel1994iterative}, \cite{rodrigues2006calculating}. Though, to our knowledge, these methods were investigated only for cases when forward problem is linear and seem to be hardly applicable for complex flow modelling. A novel technique, Adaptive Strategies PCA, was developed to address these challenges. Since linear PCA represents fields in term of truncated eigen-decomposition of its covariance matrix, AS-PCA offers two ways to adaptively adjust this representation during optimization process. This can be done either by performing rotation of the truncated eigenvectors of covariance matrix or by swaping them by more sensitive ones. To calculate optimal update, gradients of objective functions are required and using of adjoint-procedures within gradient-based optimization is proposed. Under these conditions, the method provides fast and effective reparametrization with relatively low computational cost. AS-PCA can be extended for use with other state of the art PCA-based techniques, though this question is not considered in the paper.

This paper proceed as follows. After describing the basics of gradient-based history matching in Section 2.1 and short introduction into PCA parametrization in Section 2.2, we provide theoretical and practical description of AS-PCA in Section 2.3. In Section 3, detailed results for method performance on synthetic problem are presented. We conclude with a summary and highlights for future research.

%% file: Section2/subsection2_1.tex
\subsection{History matching as inverse problem}
The inverse modelling as discussed above require finding the set of unknown model parameters \textbf{m}. In history matching, the \textbf{m} typically represented by permeability in grid blocks of reservoir model. Let $\mathbf{d_{obs}}$ represent the production data, e.g. pressure and flow rates. Then, the history matching problem can be written as follows:

\begin{equation}
\label{IPO1}
\mathbf{m^*} = \argmin_{\mathbf{m}} \big[S(\mathbf{m}) = (g(\mathbf{m}) - \mathbf{d_{obs}})^T \mathbf{C}_D^{-1} (g(\mathbf{m}) - \mathbf{d_{obs}})\big]
\end{equation}

where $(g(\mathbf{m})$ represent forward flow simulation, $\mathbf{C}_D^{-1}$ is the covariance of the observation errors.

For real reservoirs the model \textbf{m} is typically large-scale (about $10^6-10^9$ grid cells) and exceeds amount of observed data. The parameter-estimation problem is ill-posed and results in non-unique estimates, often geologically unrealistic. A common approach to make problem a well-posed one is to regularize objective function $S(\mathbf{m})$ using prior knowledge about model \cite{oliver2011recent}:

\begin{align}
\label{IPO2}
S(\mathbf{m}) = \beta (g(\mathbf{m}) - \mathbf{d_{obs}})^T& \mathbf{C}_D^{-1} (g(\mathbf{m}) - \mathbf{d_{obs}}) + \\
&+ (1 - \beta)(\mathbf{m}-\mathbf{m_{prior}})^T \mathbf{C}_M^{-1} (\mathbf{m}-\mathbf{m_{prior}}) \nonumber
\end{align}

where $\mathbf{m_{prior}}$ is prior model, $\mathbf{C}_M^{-1}$ is the covariance matrix of the prior parameter errors. $\beta$ is a scale constant.







%% file: Section2/subsection2_2.tex
\subsection{Adjoint-based history matching}
To solve the problem (\ref{IPO1}) using gradient-based methods one need an efficient procedure to calculate the gradient of objective function $S(m)$ with respect to model $m$. Because $S(m)$ consists of $g(m)$ which involves solving system of PDE, the most efficient way to obtain these gradients is to use adjoint method. In \cite{sarma2006efficient} adjoint system was defined for PDE-constrained problem. For similarity with referenced work, the problem (\ref{IPO1}) can be rewritten as:

\begin{align}
&\min_{\mathbf{m}} \Big[ S(\mathbf{m}) = \sum_{n=0}^{N-1} L^n (\mathbf{x}^{n+1}, \mathbf{m}) \Big] \qquad \forall n \in (0, ..., N-1) \\
&\text{subject to:} \nonumber \\
&g^n(\mathbf{x}^{n+1}, \mathbf{x}^n, \mathbf{m}) = 0, \qquad \forall n \in (0, ..., N-1) \nonumber\\
\end{align}

where $\mathbf{x}^n$ is the PDE solution at n-th timestep, $L^n (\mathbf{x}^{n+1}, \mathbf{m})$ is the same least squares misfit as in (\ref{IPO1}).

Following \cite{sarma2006efficient}, adjoint model can be derived as:

\begin{align}
& \mathbf{\lambda}^n = - \Big[\frac{\partial L^{n-1}}{\partial \mathbf{x}^n}+ \mathbf{\lambda}^{n+1} \frac{\partial g^n}{\partial \mathbf{x}^n} \Big] \Big[\frac{\partial g^{n-1}}{\partial \mathbf{x}^n} \Big]^{-1} \qquad \forall n \in (0, ..., N-1) \\
& \mathbf{\lambda}^N = - \Big[\frac{\partial L^{N-1}}{\partial \mathbf{x}^N} \Big] \Big[\frac{\partial g^{N-1}}{\partial \mathbf{x}^N} \Big]^{-1}
\end{align}

And finally, gradient of the objective w.r.t. the model can be calculated as:
\begin{equation}
\label{Adjoint}
\frac{dS}{d\mathbf{m}} = \sum_{n=0}^{N-1}\Big[ \frac{\partial L^{n}}{\partial \mathbf{m}} + \mathbf{\lambda}^{n+1} \frac{\partial g^{n}}{\partial \mathbf{m}} \Big]
\end{equation}

Then it can be used with any gradients-based optimization algorithm. In this paper, a nonlinear conjugate gradient (CG) method was applied.

%% file: Section2/subsection2_3.tex
\subsection{PCA-based model parametrization}
The models, generated in history matching process, have to preserve geological features of formations, i.e. to be geologically consistent. Also, performing history matching with the large scale reservoir models is a high cost process. Principal component analysis parametrization technique was developed to adress this problem. It allows to efficiently reduce model size while preserving key model features known from prior information. PCA utilize idea of stochastic process representation as an finite linear combination of orthogonal functions.

It is well known that any stochastic process can be represented as linear combination of infinite number of orthogonal functions:

\begin{equation}
\mathbf{m} = \sum_{k=1}^{\infty} A_k \boldsymbol{\varphi_k}
\end{equation}

where $\boldsymbol{\varphi_k}$ is k-th orthogonal function, $A_k$ - k-th coefficient of decomposition

The idea of data decomposition is to define such type of orthogonal basis functions, that model $\hat{m}$ given by the truncated basis approximate initial model with minimal error.

\begin{equation}
\mathbf{\tilde{m}} = \sum_{k=1}^{N} A_k \boldsymbol{\varphi_k}
\end{equation}

PCA or Karhunen-Lo\`{e}ve expansion yields the best such basis in the sense that it minimizes the total mean squared error:

\begin{equation}
\label{KLT}
    \min_{\boldsymbol{\varphi_k}} \textbf{E} \int_r  (\delta \mathbf{m(r)})^2 d\mathbf{r}
\end{equation}

where $\delta \mathbf{m}$ is defined as approximation misfit:
\begin{equation}
\label{deltam}
    \delta \mathbf{m(r)} = \mathbf{m(r)} - \mathbf{\tilde{m}(r)} =  \sum_{k=1}^{\infty} A_k \boldsymbol{\varphi_k(r)} -  \sum_{k=1}^{N} A_k \boldsymbol{\varphi_k(r)} = \sum_{k=N+1}^{\infty} A_k\boldsymbol{\varphi_k(r)}
\end{equation}


Using the fact that

\begin{equation}
    A_k = {\int_{\boldsymbol{r}} \mathbf{m(r)} \boldsymbol{\varphi_k(r)} d\mathbf{r}}
\end{equation}

Total mean square error from (\ref{KLT}) can be expressed as:
\begin{align}
\label{eq5}
    \textbf{E} \int_r  (\delta \mathbf{m(r)})^2 d\mathbf{r} = &\int_r \mathbf{E}\Bigg[ \sum_{i=N+1}^{\infty} \sum_{j=N+1}^{\infty} A_i A_j \mathbf{\varphi_i(r) \varphi_j(r)} \Bigg] d   \mathbf{r}  \nonumber \\
    =& \int_r d \mathbf{r} \sum_{i=N+1}^{\infty} \sum_{j=N+1}^{\infty} \mathbf{E}\Bigg[ \int_{s} \int_{q} \mathbf{m(s) m(q) \varphi_i(s) \varphi_j(q) ds dq} \Bigg] \mathbf{\varphi_i(s) \varphi_j(s)}  \nonumber \\
    =& \int_r d \mathbf{r} \sum_{i=N+1}^{\infty} \sum_{j=N+1}^{\infty} \mathbf{\varphi_i(s) \varphi_j(s)}  \int_{s} \int_{q} \mathbf{K(s,q) \varphi_i(s) \varphi_j(q) ds dq}  \nonumber \\
    =& \sum_{k=N+1}^{\infty} \int_s \int_q \mathbf{K(s,q) \varphi_i(s) \varphi_j(q) ds dq}
\end{align}

where $K(q,m)$ - is a covariance function defined as:
\begin{equation}
    \mathbf{K(q,m)} = \textbf{E} [\mathbf{m(s), m(q)}]
\end{equation}

Optimal orthogonal basis vector will be given by solving the following problem:
\begin{equation}
\label{KLT2}
    \int_q \mathbf{K(s,q) \varphi_i(q) dq} = \beta_i \boldsymbol{\varphi_i(s)}
\end{equation}

Equation (\ref{KLT2}) states that optimal choice for basis vectors is eigenvectors of prior models covariance matrix.

In practice we are interested in generating new models with the structure of initial model. Using KL expansion it can be performed in following manner:

\begin{equation}
\label{PCA}
    \mathbf{m} \rightarrow \mathbf{m (\xi)} =  \mathbf{\bar{m}} + \sum_{k=1}^{N} \boldsymbol{\varphi_k(x) \xi_i} = \mathbf{\bar{m}} + W \Sigma \boldsymbol{\xi}
\end{equation}

where $\mathbf{\bar{m}}$ - dataset mean model, $\boldsymbol{\xi}$ - vector from $N(0,1)$, $W$ - transition matrix with eigenvectors in columns, $\Sigma$ - diagonal matrix containing eigenvalues

Practically, $W$ and $\Sigma$ can be defined from SVD decomposition of dataset covariance matrix:

\begin{equation}
    W \Sigma V^* = \text{SVD}(K)
\end{equation}

where $K$ is a covariance matrix

Number of truncated eigenvectors is determined by energy retained in them:

\begin{equation}
\label{eq:energy}
    \text{Energy}_i = \frac{\beta_i}{\sum_i \beta_i}
\end{equation}

where $\beta_i$ is eigenvalue of corresponding eigenvector

\begin{figure}[H]
    \centering
    \includegraphics[width=\textwidth]{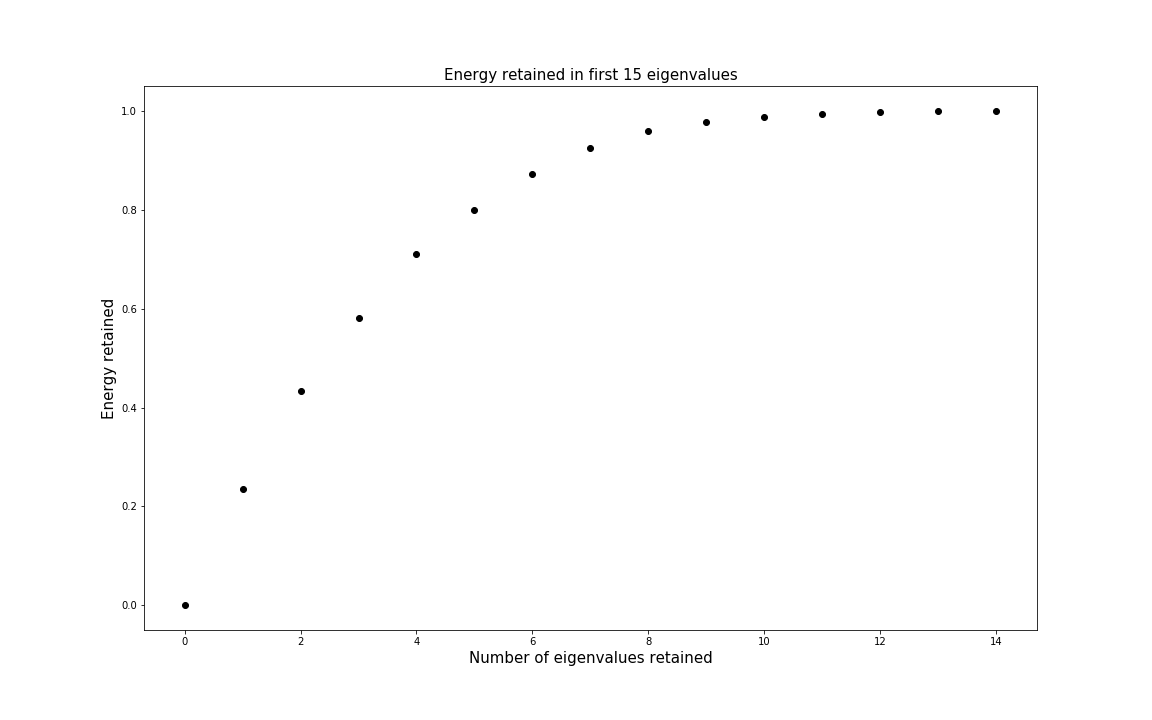}
    \caption{\textit{Energy retained in first 15 Eigenvalues}}
    \label{fig:energy}
\end{figure}

Structure of parametrization by PCA given by (\ref{PCA}) allows easily get gradient of objective function with respect to coefficients of decomposition. Combining (\ref{Adjoint}) and (\ref{PCA}):

\begin{equation}
    \frac{dS}{d\boldsymbol{\xi}} = \frac{dS}{d\mathbf{m}} \frac{d\mathbf{m}}{d\boldsymbol{\xi}} =  \Bigg[ \sum_{n=0}^{N-1}\Big[ \frac{\partial L^{n}}{\partial \mathbf{m}} + \boldsymbol{\lambda}^{n+1} \frac{\partial g^{n}}{\partial \mathbf{m}} \Big] \Bigg] (W\Sigma)^T
\end{equation}

%% file: Section3/subsection3_1.tex
\subsection{Construction of a new eigenvalue problem}
The described above principal component analysis parametrization doesn't consider sensitivity of objective function w.r.t  basis components. Therefore choice of including eigenvectors into optimization is limited by those with highest retained energy (\ref{eq:energy}). This fact leads to an incorrect model parametrization when important geological structure is underrepresented in prior models. Due to high uncertainty of data and inner reservoir structure, the problem can occur in most of real cases. In this section we propose a new technique which allows to easily include information about sensitivity into basis construction process. Our Adaptive Strategies Principal Component Analysis (AS-PCA) is much faster than sensitivity analysis process, since it uses only a few linear transforms on initial basis and doesn't require any additional run of forward model.

The basis of classic PCA is given as truncated Karhunen-Loeve expansion which provides the best approximation of the original process in the sense that it reduces the total mean-square error resulting of its truncation as described above (\ref{KLT}), (\ref{KLT2}).

Our suggestion is to find a new basis that  will minimize the misfit between objective function value of original and approximated model. 

\begin{equation}
\label{eq::minefx}
\min_{\varphi_k} \textbf{E} \Big[(\delta S)^2\Big]
\end{equation}

\begin{align}
\label{eq:defF}
\delta S(m) &\approx S(m) - S(\tilde{m}) =\nonumber \\
= S(\tilde{m} + \delta m) - S(\tilde{m})&= S(\tilde{m}) + \frac{\partial S}{\partial m(r)} \vert_{\tilde{m}} \delta m - S(\tilde{m})
\end{align}


From (\ref{eq:defF}):
\begin{equation}
\label{eq:defdf}
\delta S = \frac{\partial S}{\partial m(r)} \vert_{\tilde{m}} \delta m 
\end{equation}

Substituting definition of $\delta x$ from (\ref{deltam}):
\begin{equation}
\label{eq:df1}
\delta S = \sum_{k=N+1}^{\infty} A_k \int_r \frac{\partial S}{\partial m(r)} \varphi_k (r) dr
\end{equation}

For the sake of convenience let's define $\frac{\partial S}{\partial m(r)}$ as $J(r)$

Let's then calculate total mean square error:

\begin{equation}
(\delta S)^2 = \sum_{i=N+1}^{\infty} \sum_{j=N+1}^{\infty} A_i(\omega) A_j(\omega) \int_{r_1}\int_{r_2} J(r_1) J(r_2) \varphi_i(r_1) \varphi_j(r_2) dr_1 dr_2
\end{equation}

Using the same workflow as for classic KL expansion:

\begin{align}
(\delta S)^2 = \sum_{i=N+1}^{\infty} \sum_{j=N+1}^{\infty}& \int_s \int_q m(s) m(q) \varphi_i(s) \varphi_j(q) ds dq \times \\
& \times \int_{r_1}\int_{r_2} J(r_1) J(r_2) \varphi_i(r_1) \varphi_j(r_2) dr_1 dr_2 \nonumber
\end{align}

\begin{align}
\textbf{E}\big[(\delta S)^2 \big]= \sum_{i=N+1}^{\infty} \sum_{j=N+1}^{\infty} \int_s &\int_q K(s,q) \varphi_i(s) \varphi_j(q) ds dq \times\\ &\times\int_{r_1}\int_{r_2} J(r_1) J(r_2) \varphi_i(r_1) \varphi_j(r_2) dr_1 dr_2 \label{eq13} \nonumber
\end{align}

Then, because gradient $J(r)$ has the same shape as $r$, we can decompose it by the same basis:
\begin{equation}
\label{eq:Jbas}
J(r) = \sum_{k=1}^{\infty} c_k \varphi_k(r) \rightarrow c_k = \int_r J(r) \varphi_k(r) dr
\end{equation}

Therefore, we can define minimization problem as:
\begin{align}
\label{eq:minS}
&\min_{\varphi} \sum_{i,j} c_i c_j \int_s \int_q K(s,q) \varphi_i(s) \varphi_j(q) ds dq \\
&\text{s.t.} \int_s \varphi_i(s) \varphi_j(s) ds = 1 \nonumber
\end{align}

To find optimum basis $\varphi_k$ we will define augmented with Lagrange multipliers $\beta_i$ function:

\begin{align}
\label{eq:LA}
L_A(q) = \sum_{i,j} c_i c_j \int_s \int_q K(s,q) &\varphi_i(s) \varphi_j(q) ds dq - \\
& - \beta_i \Big[ \int_s \varphi_i(s) \varphi_j(s) ds - 1 \Big] \nonumber
\end{align}

\begin{align}
= \sum_{i} \Big[ c_i^2 \int_s \int_q K(s,q) & \varphi_i(s) \varphi_i(q) ds dq - \beta_i ( \int_s \varphi_i(s) \varphi_j(s) ds - 1) \nonumber \\
& + \sum_{j \neq i} c_i c_j \int_s \int_q K(s,q) \varphi_i(s) \varphi_j(q) ds dq \Big]
\end{align}

In order to find optimal $\varphi_k$ we need to calculate $\frac{\partial L_A}{\partial \varphi_k}$ and set it to zero:

\begin{align}
\frac{\partial L_A(r)}{\partial \varphi_i(s)} = c_i^2 \int_s \int_q K(s,q) & \varphi_i(q) ds dq - \beta_i \int_s \varphi_i(s) ds + \\
& + \sum_{j \neq i} c_i c_j \int_s \int_q K(s,q) \varphi_j(q) ds dq = 0 \nonumber
\end{align}

\begin{align}
c_i^2 \int_s \int_q K(s,q) & \varphi_i(q) ds dq +\sum_{j \neq i} c_i c_j \int_s \int_q K(s,q)\varphi_j(q) ds dq = \nonumber \\
& =\beta_i \int_s \varphi_i(s) ds 
\end{align}

Finally, we have:
\begin{align}
c_i^2 \int_q K(s,q)\varphi_i(q)dq + \sum_{j \neq i} c_i c_j\int_q K(s,q)\varphi_j(q)dq =\beta_i \varphi_i(s)
\end{align}

Dividing whole equation by $c_i^2$ and substituting $\frac{c_j}{c_i}$ by $\alpha_{ij}$ and $\frac{\beta_i}{c_i}$ by $\beta_i$ we get:
\begin{align}
\label{eq:bfQ}
\int_q K(s,q)\varphi_i(q)dq + \sum_{j \neq i} \alpha_{ij} \int_q K(s,q)\varphi_j(q)dq =\beta_i \varphi_i(s)
\end{align}

Without second addend in left side, equation (\ref{eq:bfQ}) is just the same eigenproblem as classical PCA gives (\ref{KLT2}) and solution is already known. However, the key difference here is the presence of relatively small addend, form of which doesn't allow us to easily solve the new eigenvalue problem.

%% file: Section3/subsection3_2.tex
\subsection{Stationary perturbation theory}

One can notice that equation (\ref{eq:bfQ}) looks similar to quantum mechanics problem of solving Schrodinger equation, with complex Hamiltonian $\hat{H}$ which doesn't allow to directly solve problem. Classic approach here is to divide Hamiltonian $\hat{H} = \hat{H}_0 + \hat{V}$, in such ways that eigenvectors $\varphi_k^0$ and eigenvalues $E^0$ of Hamiltonian $\hat{H}_0$ are known, and $\hat{V}$ is a small corrective of $\hat{H}_0$ or as said, perturbation and then apply Stationary Perturbation Theory (SPT) approach.

We can rewrite eigenvalue problem (\ref{eq:bfQ}) in notation that is commonly used in quantum mechanics:

\begin{equation}
\label{eq:qm1}
(\hat{H}_0 + \hat{V}) \ket{\varphi_k} = \beta \ket{\varphi_k}
\end{equation}

Known eigenvalues and eigenvectors for classic PCA case we will define as $\beta^0$ and $\varphi^0_k$.

As SPT suggests, let's decompose solution of (\ref{eq:qm1}) by the basis of standard PCA:
\begin{equation}
\label{eq:dec}
\varphi_i = \sum_k c_k^i \varphi_k^0
\end{equation}

Then, we can rewrite equation (\ref{eq:bfQ}) as:
\begin{equation}
\label{eq:SPTp}
\sum_k c_k^i \beta_k^0 \varphi_k^0 + \sum_{j \neq i} \alpha_{ij} \sum_k c_k^j \beta_k^0 \varphi_k^0 = \beta_i \sum_k c_k^i \varphi_k^0
\end{equation}

Let's multiply (\ref{eq:SPTp}) by $\varphi_n^0$ and integrate by $ds$. Then we can use that
$$\int_s \varphi_n(s) \varphi_k(s) ds = I(n=k)$$
And get:
\begin{equation}
c_n^i \beta_n^0 + \sum_{j \neq i} \alpha_{ij} c_n^j \beta_n^0 = \beta_i c_n^i
\end{equation}

Now we're going to use following representation:
\begin{align}
& \beta_i = \beta_i^0 + \beta_i^1 \\
& c_{i,k} = c_{i,k}^0 + c_{i,k}^1 \nonumber
\end{align}

where $c_k^i \rightarrow c_{i.k}$ for the sake of convenience.

From which we get:
\begin{equation}
(c_{i,n}^0 + c_{i,n}^1)\beta_n^0 + \sum_{j \neq i} \alpha_{ij} (c_{j,n}^0 + c_{j,n}^1)\beta_n^0 = (\beta_i^0 + \beta_i^1)(c_{i,n}^0 + c_{i,n}^1)
\end{equation}

Because $c_{i,n}^0 = \delta_{i,n}$, then:
\begin{align}
& i = n: \beta_n^0 + c_{n,n}^1 \beta_n^0 + \beta_n^0 \sum_{j \neq i} \alpha_{ij} c_{j,n}^0 = \beta_n^0 + \beta_n^1 + \beta_n^0 c_{n,n}^1 \Rightarrow \beta_n^1 = 0 \\
& i \neq n: c_{i,n}^0 \beta_n^0 + c_{i,n}^1 \beta_n^0 + \beta_n^0 \sum_{j \neq i} \alpha_{ij} \delta_{j,n} = \beta_i^0 c_{i,n}^0 + \beta_i^1 c_{i,n}^0 + \beta_i^0c_{i,n}^1
\end{align}

Finally, we get coefficients of decomposition (\ref{eq:dec}):
\begin{equation}
\label{eq:cj}
c_{i,n}^1 = \frac{\beta_n^0 \sum_{j \neq i} \alpha_{ij} \delta_{j,n}}{\beta_i^0 - \beta_n^0} = \frac{\alpha_{i,n} \beta_n^0}{\beta_i^0 - \beta_n^0}
\end{equation}

Eigenvectors of the problem (\ref{eq:bfQ}) are given by:
\begin{equation}
\label{eq:newbasis}
\varphi_i = \varphi_i^0 + \gamma \sum_{n \neq i} \frac{\alpha_{i,n} \beta_n^0}{\beta_i^0 - \beta_n^0} \varphi_n^0
\end{equation}

Where $\gamma$ is manually setted coefficient that will keep $V$ small relatively to $H$, as requirement for usage of SPT.
\begin{equation}
\gamma = \frac{\varepsilon}{\max \big( \sum_{n \neq i} \frac{\alpha_{i,n} \beta_n^0}{\beta_i^0 - \beta_n^0} \big)}, \qquad \varepsilon \in (0,1)
\end{equation}

%% file: Section3/subsection3_3.tex
\subsection{AS-PCA modifications}
A new parametrization basis (\ref{eq:newbasis}) is given by performing small rotation of the initial basis from classic KL expansion by using information of truncated components. Orthogonality is saved due to the smallness of addend and it allows us to efficiently embed this method into existing PCA + gradient-based optimization code just by substituting old basis by new without significant changes in workflow. Due to the fact that technique performs rotation of original PCA basis, the method was called rotation strategy of AS-PCA. 


\begin{algorithm}[H]
\KwData{PCA basis $W$, truncated PCA components $W^*$, eigenvalues $\beta$, gradient $J(r)$}
\KwResult{A new basis $\tilde{W}$}
\ForAll{$\varphi_i$ in $W$}{
calculate $c_i = \int_r J(r) \varphi_i(r) dr$ \;
\ForAll{$\varphi_j$ in $W^*$}{
calculate $c_j = \int_r J(r) \varphi_j(r) dr$\;
calculate $\alpha_{i,j} = c_i \cdot c_j$\;
calculate coefficient $c_{i,j}^1 = \alpha_{i,j} \beta_n^0 / (\beta_i^0 - \beta_n^0)$\;
update $\varphi_i = \varphi_i + c_{i,j}^1 \cdot \varphi_j$\; 
}
}
\caption{AS-PCA (rotation)}\label{rotation}
\end{algorithm}

We find interest in two heuristic modifications of idea explained above. Having total input of each (unused in the old basis) component, we can sort them by impact and then do one of two possible variants:
\begin{enumerate}
\item Extend old basis by adding some number of vectors with the highest impact, or
\item Substitute weak component of basis by them
\end{enumerate}

The first technique is called extension strategy of AS-PCA

\begin{algorithm}[H]
\KwData{PCA basis $W$, truncated PCA components $W^*$, eigenvalues $\beta$, gradient $J(r)$}
\KwResult{A new basis $\tilde{W}$}
\ForAll{$\varphi_i$ in $W$}{
calculate $c_i = \int_r J(r) \varphi_i(r) dr$ \;
\ForAll{$\varphi_j$ in $W^*$}{
calculate $c_j = \int_r J(r) \varphi_j(r) dr$\;
calculate $\alpha_{i,j} = c_i \cdot c_j$\;
calculate coefficient $c_{i,j}^1 = \alpha_{i,j} \beta_n^0 / (\beta_i^0 - \beta_n^0)$\;
save update $\varphi^*_i = c_{i,j}^1 \cdot \varphi_j$\; 
}
}
sort $W^*$ by $c_j$ in descending order\;
extend $W$ by first $n$ elements of $W^*$
\caption{AS-PCA (extension)}\label{swap}
\end{algorithm}

The procedure of substituting vectors has been named 'swap' and the whole technique is called swap strategy of  AS-PCA
\newline
\begin{algorithm}[H]
\KwData{PCA basis $W$, truncated PCA components $W^*$, eigenvalues $\beta$, gradient $J(r)$}
\KwResult{A new basis $\tilde{W}$}
\ForAll{$\varphi_i$ in $W$}{
calculate $c_i = \int_r J(r) \varphi_i(r) dr$ \;
\ForAll{$\varphi_j$ in $W^*$}{
calculate $c_j = \int_r J(r) \varphi_j(r) dr$\;
calculate $\alpha_{i,j} = c_i \cdot c_j$\;
calculate coefficient $c_{i,j}^1 = \alpha_{i,j} \beta_n^0 / (\beta_i^0 - \beta_n^0)$\;
save update $\varphi^*_i = c_{i,j}^1 \cdot \varphi_j$\; 
}
}
sort $W$ by updates $\varphi^*_i$ in descending order\;
sort $W^*$ by $c_j$ in descending order\;
swap last $n$ elements of $W$ by first $n$ elements of $W^*$
\caption{AS-PCA (swap)}\label{swap}
\end{algorithm}

%% file: Results/TPD.tex
In this section the comparative analysis of the classic PCA and proposed method is presented. A 1D non-linear diffusion problem with second-order boundary conditions was chosen as synthetic problem to test the methods:

\begin{align}
\label{difeq}
    & \frac{\partial u}{\partial t} = \frac{\partial}{\partial x} \Big( D(x) u^2 \frac{\partial u}{\partial x} \Big) \nonumber \\
    & u(x,0) = u_0, \qquad \forall x \in [0, L]  \\
    & \frac{\partial u}{\partial x}(0,t) = u_L, \qquad \forall t \in [0, T] \nonumber \\
    & \frac{\partial u}{\partial x}(L,t) = u_R, \qquad \forall t \in [0, T] \nonumber
\end{align}

By previously used nomenclature:
\begin{itemize}
    \item As a model $m$ - spatial distribution of diffusion coefficient $D(x)$;
    \item As observed data $d_{obs}$ - concentration fields $u_{obs}$, given by solving (\ref{difeq}) with known   $\hat{D}(x)$
\end{itemize}

Model size is 100 grid cells. Our goal is to define $\hat{D}(x)$ having only $u_{obs}$
\newline
\newline
The problem (\ref{difeq}) was discretized by fully implicit scheme using finite volumes method and then solved using Newton-Raphson method.

A prior information was represented by dataset of possible model realizations. Firstly, the true model was defined as:
\begin{equation}
\label{eq:true_model}
    \hat{D}(x) = 3.5 - 1.6 \sin{x} + 0.1\cos{\sqrt{300x}}
\end{equation}

Then, dataset was generated as true realization with random smooth perturbations:
\begin{figure}[H]
\centering
\includegraphics[width=9cm]{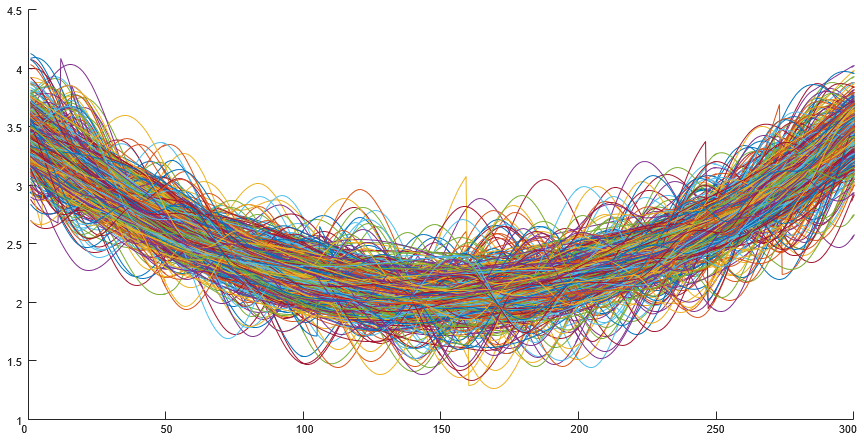}
\caption{\textit{Dataset of prior models (600 realizations))}}
\label{fig:Main}
\end{figure}

The AS-PCA algorithm was tested on two cases. In first, the true model was represented by (\ref{eq:true_model}). The second case was constructed by adding low-frequency noise into (\ref{eq:true_model})

The solutions of the problems are presented in figure below. All three strategies of AS-PCA were compared with PCA parametrization and with full model opimization without reduction. Since every basis adaptation provides new solution, we presented all of them in plots corresponding to methods.

\begin{figure}[H]
\centering
\includegraphics[width=\textwidth]{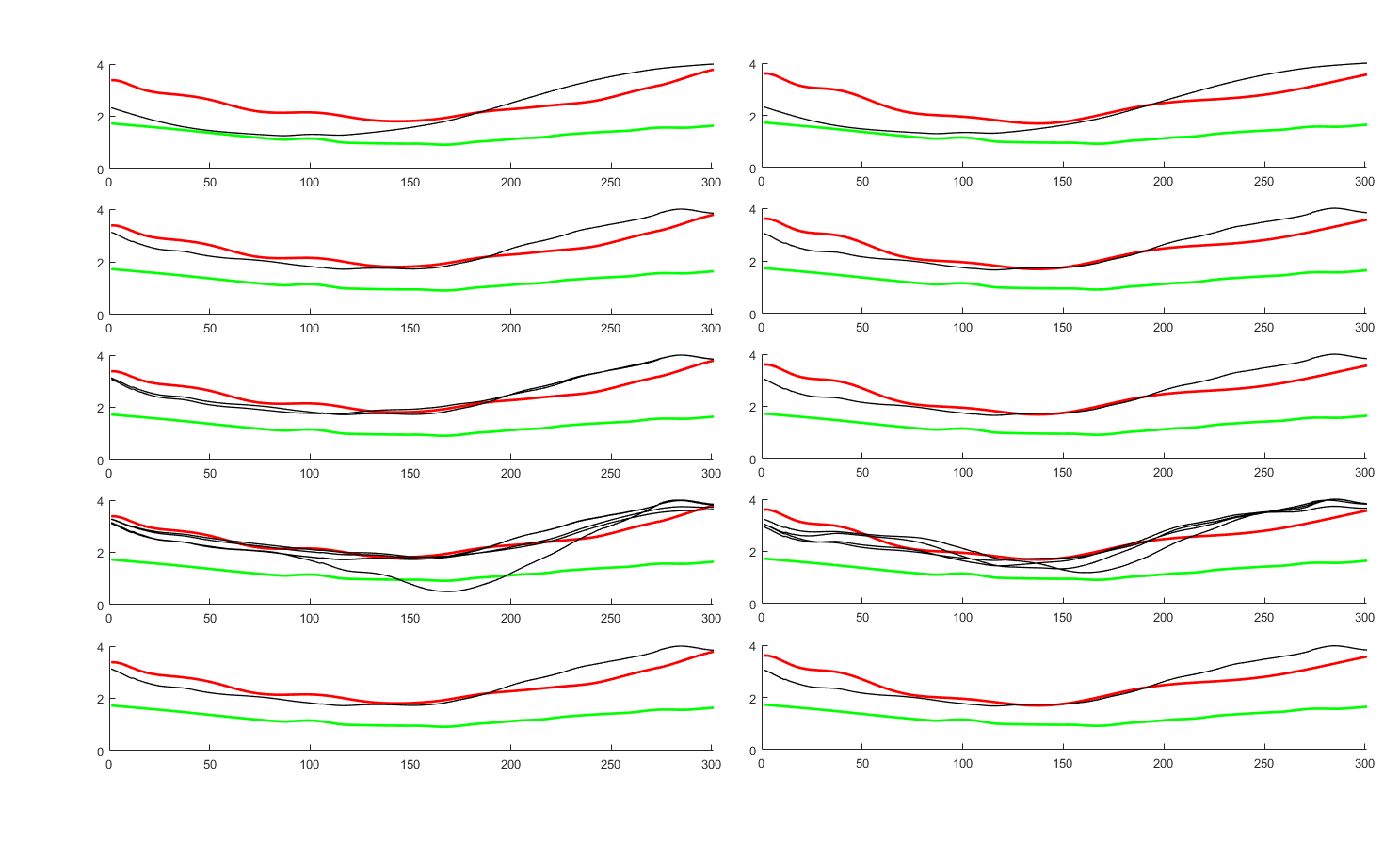}
\caption{\textit{Solutions to a synthetic problem, obtained using full size model, linear PCA, rotation, swap and extension strategies of AS-PCA: First column represents solution for true test model (\ref{eq:true_model}), second - noised test model. On all plots bold \textcolor{red}{red} line represents true model, \textcolor{forgreen}{green} - initial guess, \textbf{black} - resulting model. First row - optimization without model reduction, second - PCA parametrization, third - rotation strategy, fourth - swap strategy, fifth - extension strategy}}
\label{fig:Main}
\end{figure}

\begin{figure}[H]
\centering
\includegraphics[width=\textwidth]{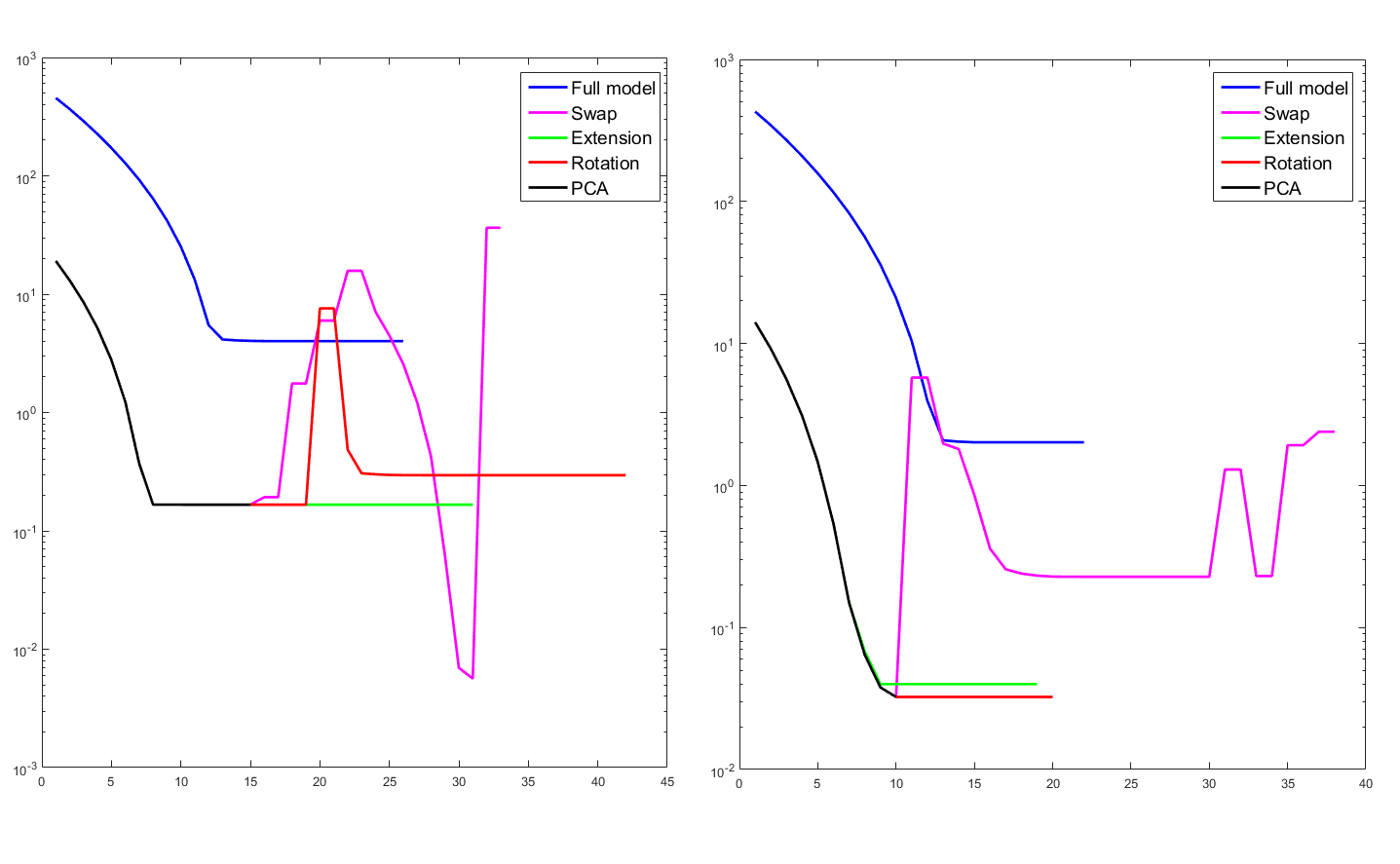}
\caption{\textit{Convergence analysis of optimization process with discussed parametrizaton algorithms. On the left - analysis for first model, on the right - for second. \textcolor{blue}{Blue} line represents convergence of objective function without parametrization, \textcolor{forgreen}{green} - extension strategy, \textcolor{red}{red} - rotation strategy, \textcolor{magenta}{magenta} - swap strategy, \textbf{black} - linear PCA}}
\label{fig:Convergence}
\end{figure}

Our general conclusions from results described below are summarized as follows:
\begin{enumerate}
    \item Rotation Strategy. A small change of basis allows to efficiently  investigate area of local minima. Our proposition is to use the strategy when all structures of solution a presumed to be well-known.
    \item Swap Strategy. Because every swap significantly change form of basis, a big amount of relatively different realizations can be obtained. Despite the fact that use of strategy doesn't always leads to a lower values of objective functions, it provides a number of realization with significantly different structures. Our assumption is to use the strategy for uncertainty quantification and forecasting, though this is topic of further research.
    \item Extension Strategy. On presented synthetic problem and dataset, significant impact of extension strategy wasn't observed. Also, method increase computational cost of problem solving, due to increasing number of parameters that have to be adjusted.
\end{enumerate}

%% file: conclusion.tex
A novel algorithm for adaptive parametrization of random fields was developed. This PCA-based method is expanded to include information about objective function structure into basis construction process. This can be done using one of three suggested strategies. One of them provide significant improve in local minima investigation (\ref{fig:Main}), second allows multiple different models obtaining. Although algorithm has been tested on synthetic examples with artificially created model complexities, applicability on more practical models is expected, since real data always have uncertain structures. Also, method provides a relatively fast parametrization, since only a few linear transforms of initial basis are required altogether with already given gradients.

Although, in this paper, use of the method is limited by gradient-based optimization, its implementation simplicity (\ref{rotation}, \ref{swap}) and relatively low computational cost allow to easily improve overall quality of inverse problem solution. Algorithm development for effective usage with non-gradient methods is one of possible directions of further research. 

Despite that AS-PCA was considered on the case of history matching problem, it is applicable for a wide class of other inverse problems, such as tomography, computer vision, etc. and can be extended for usage with a variety of parametrization techniques.

%% file: main.bbl
\begin{thebibliography}{10}

\bibitem{aanonsen2009ensemble}
S.~I. Aanonsen, G.~N{\ae}vdal, D.~S. Oliver, A.~C. Reynolds, B.~Vall{\`e}s,
  et~al.
\newblock The ensemble kalman filter in reservoir engineering--a review.
\newblock {\em Spe Journal}, 14(03):393--412, 2009.

\bibitem{hajizadeh2011towards}
Y.~Hajizadeh, M.~A. Christie, V.~Demyanov, et~al.
\newblock Towards multiobjective history matching: faster convergence and
  uncertainty quantification.
\newblock In {\em SPE reservoir simulation symposium}. Society of Petroleum
  Engineers, 2011.

\bibitem{jafarpour2008history}
B.~Jafarpour and D.~B. McLaughlin.
\newblock History matching with an ensemble kalman filter and discrete cosine
  parameterization.
\newblock {\em Computational Geosciences}, 12(2):227--244, 2008.

\bibitem{kadyrova2018application}
A.~Kadyrova and A.~Khlyupin.
\newblock Application of adjoint-based optimal control to gas reservoir with a
  memory effect.
\newblock In {\em ECMOR XVI-16th European Conference on the Mathematics of Oil
  Recovery}, 2018.
\newblock
  \href{https://arxiv.org/abs/1812.11021}{https://arxiv.org/abs/1812.11021}.

\bibitem{liu2018deep}
Y.~Liu, W.~Sun, and L.~J. Durlofsky.
\newblock A deep-learning-based geological parameterization for history
  matching complex models.
\newblock {\em arXiv preprint arXiv:1807.02716}, 2018.

\bibitem{oliver2011recent}
D.~S. Oliver and Y.~Chen.
\newblock Recent progress on reservoir history matching: a review.
\newblock {\em Computational Geosciences}, 15(1):185--221, 2011.

\bibitem{rodrigues2006calculating}
J.~R.~P. Rodrigues.
\newblock Calculating derivatives for automatic history matching.
\newblock {\em Computational Geosciences}, 10(1):119--136, 2006.

\bibitem{romero2000modified}
C.~Romero, J.~Carter, A.~Gringarten, R.~Zimmerman, et~al.
\newblock A modified genetic algorithm for reservoir characterisation.
\newblock In {\em International Oil and Gas Conference and Exhibition in
  China}. Society of Petroleum Engineers, 2000.

\bibitem{sahni2005multiresolution}
I.~Sahni, R.~N. Horne, et~al.
\newblock Multiresolution wavelet analysis for improved reservoir description.
\newblock {\em SPE Reservoir Evaluation \& Engineering}, 8(01):53--69, 2005.

\bibitem{sarma2008kernel}
P.~Sarma, L.~J. Durlofsky, and K.~Aziz.
\newblock Kernel principal component analysis for efficient, differentiable
  parameterization of multipoint geostatistics.
\newblock {\em Mathematical Geosciences}, 40(1):3--32, 2008.

\bibitem{sarma2006efficient}
P.~Sarma, L.~J. Durlofsky, K.~Aziz, and W.~H. Chen.
\newblock Efficient real-time reservoir management using adjoint-based optimal
  control and model updating.
\newblock {\em Computational Geosciences}, 10(1):3--36, 2006.

\bibitem{schulze2001optimization}
R.~Schulze-Riegert, J.~Axmann, O.~Haase, D.~Rian, Y.-L. You, et~al.
\newblock Optimization methods for history matching of complex reservoirs.
\newblock In {\em SPE Reservoir Simulation Symposium}. Society of Petroleum
  Engineers, 2001.

\bibitem{vo2014new}
H.~X. Vo and L.~J. Durlofsky.
\newblock A new differentiable parameterization based on principal component
  analysis for the low-dimensional representation of complex geological models.
\newblock {\em Mathematical Geosciences}, 46(7):775--813, 2014.

\bibitem{vo2016regularized}
H.~X. Vo and L.~J. Durlofsky.
\newblock Regularized kernel pca for the efficient parameterization of complex
  geological models.
\newblock {\em Journal of Computational Physics}, 322:859--881, 2016.

\bibitem{vogel1994iterative}
C.~R. Vogel and J.~Wade.
\newblock Iterative svd-based methods for ill-posed problems.
\newblock {\em SIAM Journal on Scientific Computing}, 15(3):736--754, 1994.

\end{thebibliography}
